\documentclass{article}
\usepackage[total={6in, 10in}]{geometry}
\usepackage{amssymb,amsmath,bm}
\usepackage{amsthm}
\usepackage{color}
\usepackage{mathrsfs,amsmath}
\usepackage[version=4]{mhchem}
\usepackage{mathrsfs}
\usepackage{graphicx}
\usepackage{bm}
\usepackage{amsfonts}
\usepackage{multirow}
\usepackage[ruled]{algorithm2e}
\usepackage{amssymb,amsmath,bm}
\usepackage{amsmath}
\usepackage{amssymb}
\usepackage{amsthm}
\usepackage{longtable}
\usepackage{multirow}
\usepackage{setspace}
\usepackage{booktabs}
\usepackage{color}
\usepackage{mathrsfs,amsmath}
\usepackage[version=4]{mhchem}
\usepackage{mathrsfs}
\usepackage{graphicx}
\usepackage{dcolumn}
\usepackage{bm}
\usepackage{natbib}
\usepackage{amsthm}
\usepackage{ulem}
\usepackage{setspace}
\usepackage{color}
\usepackage{amsfonts}
\usepackage{amsmath}
\usepackage{multirow}
\usepackage{hhline}  
\usepackage{booktabs}
\usepackage[breaklinks=true,colorlinks=true,linkcolor=blue,urlcolor=blue,citecolor=blue]{hyperref}
\usepackage{subfig}
\usepackage{float}  
\usepackage{authblk}
\usepackage{natbib}
\setcitestyle{numbers,square}
\newtheorem{thm}{Theorem}[section]

\newtheorem{rem}[thm]{Remark}

\providecommand{\keywords}[1]{\textbf{Keywords:} #1}

\title{Entropy Structure Informed Learning for Inverse XDE Problems}
\date{} 
\author[ ]{Yan Jiang\textsuperscript{2$,*$},Wuyue Yang$^{1,}$\thanks{These authors have contributed equally to this work}$^,$\thanks{yangwuyue@bimsa.cn}}
\author[ ]{Yi Zhu$^{1,2}$}
\author[ ]{Liu Hong$^{3,}$\thanks{hongliu@sysu.edu.cn}}
\affil[1]{Yanqi Lake Beijing Institute of Mathematical Sciences and Applications, Beijing 101408, China.}
\affil[2]{Yau Mathematical Sciences Center, Tsinghua University, Beijing, 100084, China.}
\affil[3]{School of Mathematics, Sun Yat-Sen University, Guangzhou, 510275, China.}

\begin{document}
\maketitle
\begin{abstract}
Entropy, since its first discovery by Ludwig Boltzmann in 1877, has been widely applied in diverse disciplines, including thermodynamics, continuum mechanics, mathematical analysis, machine learning, etc. In this paper, we propose a new method for solving the inverse XDE (ODE, PDE, SDE) problems by utilizing the entropy balance equation instead of the original differential equations. This distinguishing feature constitutes a major difference between our current method and other previous classical methods (e.g. SINDy). Despite concerns about the potential information loss during the compression procedure from the original XDEs to single entropy balance equation, various examples from MM reactions, Schlögl model and chemical Lorenz equations in the form of ODEs to nonlinear porous medium equation and Fokker-Planck equation with a double-well potential in the PDE form all well confirm the accuracy, robustness and reliability of our method, as well as its comparable performance with respect to SINDy.
\end{abstract}



\section{Introduction}
Data-driven equation discovery aims to explore an explicit closed form of equations based on time series data obtained from measurements. Traditionally, the governing equations are derived based on the first principle, like conservation laws of energy and moment, meanwhile the unknown coefficients are estimated by the linear regression method \cite{weisberg2005applied}. Nowadays, given the great power of neural networks for approximating any continuous function easily and especially in the high-dimensional analysis, machine learning algorithms have been extensively applied to data-driven equation discoveries. For example, Brunton et al. proposed a sequential threshold ridge regression algorithm (STRidge)  for revealing nonlinear dynamic systems which was an early influential work in this field \cite{Brunton2016discoverying}. Sparse regression also extends to partial differential equations (PDE-FIND) and parametric partial differential equations \cite{doi:10.1126/sciadv.1602614,rudy2019data}. Due to their easy implementation and accuracy, physics-informed neural networks have been widely used to solve the forward and inverse problems of differential equations\cite{raissi2019physics,CiCP-29-698,both2021deepmod,lu2021learning}.  Chen et al. present a novel method (PINN-SR) by combining PINN with sparse regression to achieve impressive results\cite{chen2021physics}. A far-from-complete list on machine learning based algorithms for discovering differential equations is summarized in Table \ref{MethodsSummary}.



Many studies found that a ``nice'' structure of differential equations offers a lot of help during the solvation of the inverse problems of XDEs. Among them, the general principle of thermodynamics provides insightful guidance for constructing such a kind of structures. For example, data-driven machine learning algorithms that can preserve either the Hamiltonian structure \cite{greydanus2019hamiltonian,zhong2019symplectic} or the Lagrangian structure \cite{cranmer2020lagrangian}, have been applied to a number of cases with success. Yu et al. utilize a generalized Onsager's principle to guide the automatic discovery of dissipative dynamics and their governing equations \cite{yu2021onsagernet}. With an emphasis on energy conservation and entropy dissipation based on GENERIC, a thermodynamically consistent deep neural network  \cite{hernandez2021deep} has been proposed to predict physical variables and their corresponding dynamics. Integration of neural networks and thermodynamically compatible hyperbolic equations with applications to heat conditions and one-component fluids can be found in Refs.  \cite{ZHAO2022122396,HuangMaZhouYong+2021+355+370}. 

\begin{table}[h]
\centering
\setlength\extrarowheight{10pt}
\setlength{\tabcolsep}{2pt}
\begin{tabular}{|c|c|l|}
\hline
\begin{tabular}[c]{@{}c@{}}Differential \\ equations\end{tabular} & Methods                                                                  & \multicolumn{1}{c|}{Brief description}                                                                                                                      \\ \hline
\multirow{7}{*}{ODE}                                              & SINDy\cite{Brunton2016discoverying}                     & Sparse regression and least squares techniques.                                                        \\ \cline{2-3} & DeepXDE\cite{lu2021deepxde}& An integrated python library \\ \cline{2-3} & OnsagerNet\cite{yu2021onsagernet} & \begin{tabular}[c]{@{}l@{}}Based on the generalized Onsager principle,\\ neural networks can preserve physical structure information\end{tabular} \\ \cline{2-3}& ODENet\cite{HU2022111203} & Neural ODE combines with sparse regression\\ \cline{2-3}& CRNN\cite{ji2021autonomous} & Chemical Reaction Neural Network \\ \cline{2-3}& SymODEN\cite{zhong2019symplectic}& Symplectic ODE-Net \\ \cline{2-3}  & Lift \& Learn\cite{qian2020lift}  & Identify a coordinate transformation\\ \hline
\multirow{6}{*}{PDE}                                                                & PINN\cite{raissi2019physics}                            & Physics-Informed neural network                              \\ \cline{2-3}
  & PDE-FIND\cite{doi:10.1126/sciadv.1602614}   & Sparse regression and least squares techniques. \\ \cline{2-3}
 & PINN-SR\cite{chen2021physics} & PINN and sparse regression  \\ \cline{2-3}
& SGA-PDE\cite{chen2021symbolic}  & Symbolic genetic algorithm  \\ \cline{2-3}
\multicolumn{1}{|l|}{}                                            & PDE-Net 1.0, PDE-Net 2.0\cite{long2018pde,long2019pde} & A new feed-forward deep network                                                                                                                             \\ \cline{2-3}
\multicolumn{1}{|l|}{}                                            & DLGA-PDE\cite{xu2020dlga}                               & Combination of deep learning and genetic algorithm                                                                      \\ \cline{1-3}  
\multirow{4}{*}{SDE}                                              & PINN\cite{chen2021solving}                              & Using Fokker-Planck equation and PINN                                             \\ \cline{2-3} 
  & SINDy\cite{callaham2021nonlinear} & Langevin regression   \\ \cline{2-3}  & Weak form of FPE\cite{ma2021learning}                   & Weak form of Fokker planck equation and WGAN    \\ \cline{2-3} & UQ with PINN\cite{zhang2019quantifying}                 & Using uncertainty quantification and PINN   \\ \hline
\end{tabular}
\caption{A summary of machine learning based methods for solving the forward and inverse problems of XDEs.}
\label{MethodsSummary}
\end{table}

Motivated by previous works, especially those utilized the thermodynamic principles and thermodynamic structures, here we propose a new method for solving the inverse problem of XDEs based on the entropy function and its balance equation. Once the form of an XDE is fixed (may still contain unknown parameters and unknown terms), the time evolution equation of the entropy or relative entropy could be derived with respect to the XDEs. In this way, the dynamics of XDEs is embedded into its entropy balance equation. Consequently, a natural question arises, could we extract a full knowledge of all unknown parameters in the XDEs from the time series data of entropy function in most cases? 

We aim to give a positive answer in this work. 
For a number of classical examples in mathematical biology and mathematical physics ranging from ODEs to PDEs and SDEs, we show that the time series data of entropy function indeed contains adequate information for solving the inverse problems of XDEs. An expanding list includes the Michaelis–Menten kinetics, Schlögl model, chemical Lorenz equations, nonlinear 2D porous medium equation and stochastic Langevin-type dynamics with a double-well potential. The merit of our entropy structure informed learning (EnSIL) is highlighted under different circumstances by varying the time step, noise ratio, initial conditions, etc., and compared with some of the state-of-the-art machine learning methods, such as SINDy and PDE-FIND.

The remainder of the paper is organized as follows. Section 2 establishes the theoretical foundation and applicable procedure of entropy-structure informed method. The EnSIL method is then applied to a number of classical examples, including ordinary differential equations, partial differential equations and stochastic differential equations. The details of application can be found in Section 3. The last section contains the conclusion and some discussions related to the EnSIL method.

\section{Methodology}
\subsection{Entropy and entropy balance}
Here we consider the evolutionary XDEs (ODEs and PDEs) in a general form  
\begin{eqnarray}
\label{de}
&\frac{\partial \mathbf{u}}{\partial t}+\mathcal{G}(\mathbf{x},t;\mathbf{u},\mathbf{u}_\mathbf{x}, \mathbf{u}_{\mathbf{x}\mathbf{x}},\cdots)=0,
\end{eqnarray}
where $\mathbf{u}=(u_1(\mathbf{x},t),u_2(\mathbf{x},t), \cdots, u_n(\mathbf{x},t))^T$ is an $n$-dimensional vector, denoting a solution to the above equation accompanied with proper initial/boundary conditions at time $t$ and spatial position $\mathbf{x}=(x_1,x_2,\cdots,x_d)^T$ in the  domain $t\in\mathbb{R}^1, \mathbf{x}\in\Omega \subseteq \mathbb{R}^d$. $\mathcal{G}$ is a continuous functional of 
$(\mathbf{x},t;\mathbf{u},\mathbf{u}_\mathbf{x},\mathbf{u}_{\mathbf{x}\mathbf{x}},\cdots)$ with $\mathbf{u}_\mathbf{x}$ and $\mathbf{u}_{\mathbf{x}\mathbf{x}}$ being the Jacobian and Hessian matrices.

It is well known that for many different forms of XDEs, an  entropy function $S=S(\mathbf{u})$ as a functional of the solutions $\mathbf{u}(\mathbf{x},t)$ could be introduced. This entropy function enjoys many significant merits, like the concavity with respect to variables $\mathbf{u}$, the monotonicity in the time domain, and so on. Most importantly, the time evolutionary information of $\mathbf{u}$ could be incorporated into a single balance equation for entropy, i.e.
\begin{eqnarray*}
    \frac{dS}{dt}=\frac{\partial S}{\partial \mathbf{u}}\cdot \frac{\partial \mathbf{u}}{\partial t}= - S_\mathbf{u}\cdot \mathcal{G}(\mathbf{x},t;\mathbf{u},\mathbf{u}_\mathbf{x}, \mathbf{u}_{\mathbf{x}\mathbf{x}},\cdots)
    =\frac{d^e S}{dt}+ \frac{d^i S}{dt}.
\end{eqnarray*}
The last equality represents the fact that the entropy change in time could be decomposed into rates of entropy flux $\frac{d^e S}{dt}$ and entropy production $\frac{d^i S}{dt}$ based on a knowledge of thermodynamics. Furthermore, it is remarkable that the entropy production rate (EPR) is always non-negative, $\frac{d^i S}{dt}\geq 0$, according to the \textit{Second Law of Thermodynamics}. The EPR is zero if and only if the system attains its equilibrium state.  

\subsection{Entropy balance equation for ODEs}
When the system is spatially homogeneous, we arrive at a simple version -- a group of first-order ordinary differential equations, 
\begin{eqnarray}
\label{ode}
&&\frac{d \mathbf{u}}{dt}+\mathcal{G}(t,\mathbf{u})=0,\\
&&\mathbf{u}(0)=\mathbf{u}_0.\nonumber
\end{eqnarray}
It is known that all kinds of high-order ODEs can be cast into the above form by a change of variables. 
With respect to a concrete form of Eq. \eqref{ode}, the existence of a Lyapunov function $F(t)=F(\mathbf{u}(t))$ with elegant mathematical properties stated below can often be shown. 
\begin{rem}
The Lyapunov function of Eq. \eqref{ode} has properties as follows,
\begin{itemize}
\item[(1)]$F(\mathbf{u})$ is convex;
\item[(2)]$F(\mathbf{u})\geq 0$ and $L(\mathbf{u}) = 0$ if and only if the system reaches its equilibrium state.
\item[(3)]$d F(\mathbf{u}) / d t \leq 0$.
\end{itemize}
\end{rem}
\noindent In the literature, the Lyapunov function is also known as the relative entropy.

As an example, we consider a general isolated chemical reaction system with $N$ species and $M$ reversible  reactions under constant temperature and constant volume,
\begin{equation}
\frac{dc_j}{dt}=\sum_{l=1}^M (\nu_{lj}^- - \nu_{lj}^+)\left(r_l^+(\mathbf{c}) -r_l^-(\mathbf{c})\right),\quad j=1-N,
\label{mass-action}
\end{equation}
where $\mathbf{c}(t)=(c_1(t), c_2(t), \cdots, c_N(t))^T$ denote the concentrations of chemical species at time $t$, $\nu_{lj}^+$ and $\nu_{lj}^-$ are the stoichiometric coefficients, $r_l^+(\mathbf{c})$ and $r_l^-(\mathbf{c})$ are the forward and backward chemical reaction rates for the $l^{th}$ reaction, respectively. 
Specifically, by adopting the law of mass action, the reaction rates become a form of polynomials,
\begin{equation}
r_l^+(\mathbf{c})=k_{l}^+ \prod_{i=1}^N (c_i)^{\nu_{lj}^+}, \quad r_l^-(\mathbf{c})=k_{l}^- \prod_{i=1}^N (c_i)^{\nu_{lj}^-},
\end{equation}
here $k_{l}^+$ and $k_{l}^-$ are the forward and backward rate constants.

Provided the chemical reaction system in Eq. \ref{mass-action} satisfies the condition of detailed balance (or complex balance), a Lyapunov function
\begin{equation}
F(t)=\sum_{j=1}^{N}\left(c_{j} \ln \frac{c_{j}}{c_{j}^{s s}}-c_{j}+c_{j}^{s s}\right),
\label{REF}
\end{equation} 
can be introduced, where $(c_{1}^{s s}, c_{2}^{s s}, \cdots, c_{N}^{s s})$ denote the equilibrium (or static) solution. In chemistry, the above Lyapunov function is called the Helmholtz free energy.

Accordingly, the free energy dissipation rate is given by
\begin{eqnarray}
\frac{d F(t)}{d t}
&=&\sum_{j=1}^{N}\frac{dc_{j}}{dt} \ln \left(\frac{c_{j}}{c_{j}^{s s}}\right)\nonumber\\
&=&-\sum_{l=1}^{M}\left(r_l^+(\mathbf{c}) -r_l^-(\mathbf{c})\right)\left[\ln \frac{r_l^+(\mathbf{c})}{r_l^+(\mathbf{c}^{ss})}- \ln \frac{r_l^-(\mathbf{c})}{r_l^-(\mathbf{c}^{ss})}  \right]\leq 0,
\label{relative_energy}
\end{eqnarray}
which vanishes if and only if $c(t)=c^{ss}$. This is a reminiscent of the \textit{Second Law of Thermodynamics}.

\subsection{Entropy balance equation for PDEs}
Compared to ODEs, partial differential equations are far more complicated. For simplicity, here we only illustrate the results for hyperbolic equations based on the Conservation-Dissipation Formalism (CDF). The basic model is a group of first-order hyperbolic equations, which read
\begin{equation}
\label{cdf}
\frac{\partial }{\partial t}
\left(
  \begin{array}{ccc}
    \mathbf{u}_1\\
    \mathbf{u}_2
  \end{array}
\right)      
+\sum_{l=1}^d \frac{\partial }{\partial x_l}
\left(
  \begin{array}{ccc}
    \mathbf{f}_l(\mathbf{u}_1,\mathbf{u}_2)\\
    \mathbf{g}_l(\mathbf{u}_1,\mathbf{u}_2)
  \end{array}
\right)=\left(
  \begin{array}{ccc}
    0\\
    \mathbf{q}(\mathbf{u}_1,\mathbf{u}_2)
  \end{array}
\right)
,
\end{equation}
where $(\mathbf{u}_1(\mathbf{x},t), \mathbf{u}_2(\mathbf{x},t))^T$ represent the conserved and dissipative variables at time $t$ and position $\mathbf{x}\in\Omega\subseteq\mathbb{R}^d$, $(\mathbf{f}_l, \mathbf{g}_l)^T$ are the corresponding fluxes respectively,  $\mathbf{q}(\mathbf{u}_1,\mathbf{u}_2)$ is the source term. 

To ensure the long-time asymptotic stability of the solutions, CDF suggests a sufficient condition, that is (i) the existence of a strict concave local entropy function $s(\mathbf{u}_1,\mathbf{u}_2)$ satisfying the requirement of $s_{\mathbf{u}\mathbf{u}}\cdot (\mathbf{f}_{l\mathbf{u}}, \mathbf{g}_{l\mathbf{u}})^T$ to be symmetric for each $l$; (ii) the existence of a positive definite dissipation matrix $\mathbf{M}(\mathbf{u}_1,\mathbf{u}_2)$ such that $\mathbf{q}=\mathbf{M}\cdot s_{\mathbf{u}_2}$.

According to above conditions, the time derivative of local entropy function satisfies
\begin{eqnarray*}
\frac{d s}{dt}
&&=\frac{\partial s}{\partial \mathbf{u}_1}\cdot \frac{\partial \mathbf{u}_1}{\partial t} + \frac{\partial s}{\partial \mathbf{u}_2}\cdot \frac{\partial \mathbf{u}_2}{\partial t}= s_{\mathbf{u}_1}\cdot (-\nabla_x\cdot\mathbf{f}) + s_{\mathbf{u}_2}\cdot (-\nabla_x\cdot\mathbf{g}+\mathbf{q})\\
&&=-(s_{\mathbf{u}_1}\cdot \nabla_x\cdot\mathbf{f} + s_{\mathbf{u}_2}\cdot \nabla_x\cdot\mathbf{g})+s_{\mathbf{u}_2}\cdot\mathbf{q}=-\nabla_x\cdot \mathbf{J}^s+s_{\mathbf{u}_2}\cdot\mathbf{M}\cdot s_{\mathbf{u}_2},
\end{eqnarray*}  
where $\mathbf{f}=(\mathbf{f}_1, \cdots, \mathbf{f}_d)$, $\mathbf{g}=(\mathbf{g}_1, \cdots, \mathbf{g}_d)$. The last equality holds by applying the stability condition of CDF and Poincare's Lemma. The above formula is recognized as the local entropy balance equation for hyperbolic equations in Eq. \ref{cdf}, since the first term on the right-hand side is the local entropy flux while the second term is the local entropy production rate, $s_{\mathbf{u}_2}\cdot\mathbf{M}\cdot s_{\mathbf{u}_2}\geq0$. In addition, by introducing the total entropy function $S(t)=\int_{\Omega}s(\mathbf{u}_1,\mathbf{u}_2)d\mathbf{x}$ and using the method of integration by parts, we can show the entropy balance equation holds in a global sense too.

Similar results are applied to the parabolic equations too, which will not be addressed here. Interested readers may refer to Ref. \cite{evans2004entropy}.

\subsection{Entropy balance equation for SDEs}
In the presence of either intrinsic randomness or external noises, many physical processes could be the Langevin equation, which reads
\begin{equation}
\label{le}
d\mathbf{x}(t)=\mathbf{u}(\mathbf{x})dt+\mathbf{\sigma}(\mathbf{x})d\mathbf{W}(t),
\end{equation}
where the vector $\mathbf{x}=\mathbf{x}(t) \in \mathbb{R}^n$ is the state variable. It evolves according to two different contributions -- a deterministic part $\mathbf{u}(\mathbf{x})$ and a random part $\mathbf{W}(t)$. The latter, known as the Wiener process in mathematics or the Brownian motion in physics, is a multi-dimensional stochastic process satisfying $\mathbf{W}(0)=0$ with probability $1$, having continuous paths, stationary and independent Gaussian increments.

Each solution of Eq. \eqref{le} corresponds to a single trajectory in space. An ensemble of such trajectories gives the time evolution of the probability density function (pdf) $p(\mathbf{x},t)$ rather than the stochastic variable $\mathbf{x}(t)$ itself. Based on the Itô's formula, it is straightforward to derive the governing equation of the pdf $p(\mathbf{x},t)$ -- the Fokker-Planck equation as
\begin{equation}
\label{fp}
{\partial_t}{p(\mathbf{x},t)}=-\nabla_x \cdot \mathbf{J}(\mathbf{x},t),\quad
\mathbf{J}(\mathbf{x},t)=\mathbf{u}(\mathbf{x})p(\mathbf{x},t)-\mathbf{D}(\mathbf{x})\cdot \nabla_x p(\mathbf{x},t),
\end{equation}
where $\mathbf{J}(\mathbf{x},t)$ serves as the probability flux with $\mathbf{D}(\mathbf{x})=\mathbf{\sigma}(\mathbf{x})^T\mathbf{\sigma}(\mathbf{x})/2$ representing the diffusion matrix. As to the FPE, an entropy function can be introduced, $S(t)=-\int d\mathbf{x} p(\mathbf{x}, t) \ln p(\mathbf{x}, t)$. The corresponding entropy balance equation reads
\begin{equation}
\frac{d S(t)}{d t}=-\int \mathbf{J}(\mathbf{x}, t) \nabla_x \ln p(\mathbf{x}, t) d\mathbf{x},
\end{equation}
by applying the method of integration by parts and boundary conditions at infinity.

Further suppose the FPE in \eqref{fp} has a stationary solution $p^{ss}(\mathbf{x})$, which solves the equation $\nabla \cdot J^{ss}(x)=0$ with the stationary flux $\mathbf{J}^{ss}(\mathbf{x})=\mathbf{u}(\mathbf{x})p^{ss}(\mathbf{x})-\mathbf{D}(\mathbf{x})\cdot \nabla_x p^{ss}(\mathbf{x})$. Now it becomes possible to introduce a new free energy function (or the relative entropy function)
\begin{equation}
F(t)=\int p(\mathbf{x}, t) \ln \frac{p(\mathbf{x},t)}{p^{ss}(\mathbf{x})}d\mathbf{x}\geq0,
\end{equation}
which is monotonically decreasing with time $t$, i.e.
\begin{equation}
\frac{d F(t)}{d t}=\int \mathbf{J}(\mathbf{x}, t) \nabla_x \ln \frac{p(\mathbf{x}, t)}{p^{ss}(\mathbf{x})}\leq0.
\end{equation}
Correspondingly, the entropy balance equation can be reformulated into
\begin{equation}
\frac{d S(t)}{d t}=-\int \mathbf{J}(\mathbf{x}, t) \nabla_x \ln p^{ss}(\mathbf{x}) d\mathbf{x}+\int \mathbf{J}(\mathbf{x}, t) \nabla_x \ln \frac{p^{ss}(\mathbf{x})}{p(\mathbf{x}, t)}.
\end{equation}
The two terms on the right-hand side are known as the excess heat and non-adiabatic entropy production rate in the field of stochastic thermodynamics.

\subsection{Algorithm implementation}
How to numerically calculate terms on the right-hand side of the entropy balance equation with high precision constitutes  a major bottleneck during the application of the current method. To solve this issue, we refer to the Legendre polynomials, which provide computational speed with low precision loss\cite{cohen2012polynomial} and can approximate arbitrary continuous functions within the region $[-1,1]$. As a consequence, we rewrite the entropy balance equation as
\begin{equation}\label{dFdt}
\left\{\begin{array}{l}
\displaystyle\frac{d S}{d t}=P_{\Theta}(t),\\
P_{\Theta}(t)=\sum_{m=1}^n \Theta^{(m)} \cdot \text{LPA}_m(t).
\end{array}\right.
\end{equation}
whose terms on the right-hand side are replaced by parameterized Legendre polynomials $P_{\Theta}$ truncated to the $n^{th}$  order as a function of time $t$. To be concrete, the $m^{th}$  order Legendre polynomial can be expressed as
\begin{equation}
\text{LPA}_m(t)=\frac{1}{2^m} \sum_{s=0}^{[m / 2]} \frac{(-1)^s(2 m-2 s) !}{s !(m-s) !(m-2 s) !} t^{m-2 s},\quad t \in[-1,1].
\end{equation}  
And $\Theta^{(m)}$ is a free parameter corresponding to the $m^{th}$ term. 

To reduce the impact of external noises as much as possible, we'd better rephrase Eq. \eqref{dFdt} into an integral form, i.e. $S(t)=S(0)+\int_0^t P_{\Theta}(\tau)d\tau$. See Refs. \cite{schaeffer2017sparse,kaheman2020sindy,messenger2021weak} for details. Especially, with respect to the Legendre polynomials we adopted, the integration can be carried out in an explicit way. It exhibits a dramatic advantage of our above method that there is no numerical error for integration. We only have to pay attention to the approximation errors.

As a key step, the free parameters $\Theta$ can be determined through optimizing the loss function $\mathcal{L}(\Theta)$, which includes three parts -- data loss, thermodynamic loss and constraint loss.
\begin{equation}
\mathcal{L}(\Theta)=\omega_1\mathcal{L_{\text{data}}}+\omega_2\mathcal{L_{\text{thermodynamics}}}+\omega_3\mathcal{L_{\text{constraint}}},
\label{loss}
\end{equation}
where
\begin{equation}
\begin{aligned}
& \mathcal L_{\text {data }}=\frac{1}{M} \sum_{i=1}^M\left(\hat{S}\left(t_i\right)-S\left(t_i\right)\right)^2, \\
& \mathcal L_{\text {thermodynamics}}=\bigg\|\text{ReLU}\bigg(-\frac{d^iS}{dt}\bigg)\bigg\|_2^2, \\
& \mathcal L_{\text {constraint}}=\bigg\|\bigg(\frac{dS}{dt}\bigg|_{t_N}\bigg)\bigg\|_2^2.
\end{aligned}
\end{equation}
The first loss term $\mathcal L_{\text {data }}$ is associated with the difference between the predicted entropy function $\hat{S}(t_i)=S(0)+\int_0^{t_i} P_{\Theta}(\tau)d\tau$ and real one derived from the sample data $S(t_i)$. The second loss term $\mathcal L_{\text {thermodynamics }}$ characterizes the violation of the predicted time trajectories on the second law of thermodynamics since $\frac{d^iS}{dt}\geq 0$ for the real data. The last loss term $\mathcal L_{\text {constraint}}$ will be adopted for a long-time dynamics to ensure the system approaches to its equilibrium (or static) state. $\omega_1, \omega_2, \omega_3$ are hyper-parameters.


Finally, in this paper we adopt the mean relative error (MRE) to evaluate the accuracy of predictions
\cite{doi:10.1126/sciadv.1602614}. 
\begin{equation}
\text { MRE }=\frac{1}{n} \sum_{j=1}^n\left|\frac{\Theta_j-\hat{\Theta}_j}{\Theta_j}\right| \times 100 \%, 
 \Theta_j\neq0,
\end{equation}
where $\Theta_j$ and $\hat{\Theta}_j$ represent the true and predicted parameters, respectively. $n$ is the total number of nonzero free parameters to be determined.


\subsection{Error analysis}
To gain a better understanding of the advantages and limitations of the EnSIL method, especially to make a direct comparison with classical machine learning algorithms based on original XDEs, we perform the error analysis on the Fokker-Planck equation as an illustration. With respect to the data loss $[S(t)-S(0)-\int_0^td\tau\int dx J(x,t) \nabla_x \ln p(x,t)]^2$, the errors during applying the EnSIL method come from three main sources.

\begin{itemize}
    \item \textbf{Statistical errors:} A key step in applying the EnSIL method to the Fokker-Planck equation is to estimate the probability density $p(x,t)$, which heavily relies on the sample data size and its quality. The difference between the empirical probability density $\tilde{p}(x,t)$ and the real one will cause statistical errors during the calculation of the entropy function, i.e. $$S(t)=\int p(x,t)\ln p(x,t)dx\approx \int \tilde{p}(x,t)\ln \tilde{p}(x,t)dx.$$ 
    \item \textbf{Spatial derivative errors:} Both flux term $J(x,t)$ and force term $\nabla_x \ln p(x,t)$ involve spatial derivatives. To improve the precision, the spectral method is adopted, whose errors often exhibit exponential convergence with respect to the data size. Compared to usual finite difference scheme, spectral methods can also filter noises and have better approximation properties. $$\nabla_x p(x,t)\approx\mathcal{FFT}^{-1}(iw\mathcal{FFT}(p(x,t))),$$
    where $\mathcal{FFT}$ stands for the fast Fourier transform.
    \item \textbf{Spatial integration errors:} 
    To calculate either the entropy function or its time derivative, spatial integration is unavoidable. Here we adopt the classical trapezoidal rule for the 1d case and the Monte Carlo integration for high dimensions. The scheme for the latter reads 
    \begin{eqnarray*}
    &&S(t)=\int dxp(x,t)\ln p(x,t)\approx \sum_{i=1}^N \ln p(x_i,t),\\
    &&\int dx J(x,t) \nabla_x \ln p(x,t)\approx \sum_{i=1}^N J(x_i,t) \frac{\nabla_x p(x_i,t)}{p^2(x_i,t)},
    \end{eqnarray*}
    in which the data points $x_i$ are sampled according to the distribution $p(x,t)$. And it is well-known that the error for the Monte Carlo integration decreases with the sample size in a rate of $N^{-1/2}$.
    \item \textbf{Legendre truncation errors:} To make a better estimation on the temporal integration, we adopt the Legendre polynomial approximation, which means the function below will be approximated through the $m^{th}$ order Legendre polynomials,
    \begin{equation*}
    \int dx J(x,t) \nabla_x \ln p(x,t)\approx LPA_m(t).
    \end{equation*}
    As the temporal integration of Legendre polynomials could be exactly done, we only need to consider the truncation error, which will be gradually decreasing with an increase of $m$. Meanwhile, a higher computational cost incurs.
\end{itemize}

\section{Numerical results}
In this section, we will apply our Entropy Structure Informed Learning (EnSIL) to solve the inverse problems related to ODEs, PDEs and SDEs separately. Three representative ODE models in the field of chemical reactions were tested. The first case is the Michaelis-Menten kinetics, in which we aim to explore the influence of different time steps on the results. The second case is a simple yet remarkable example of a reaction network exhibiting bi-stability -- the Schlögl model. Through this model, the robustness of EnSIL in the presence of data noises is highlighted. In the third case, we choose the chemical Lorenz equations, a complicated multi-scale open reaction system with the rate constants varying over eight orders of magnitude. Furthermore, this system is chaotic too. The nonlinear porous medium equation and the Langevin dynamics with a double-well potential have been considered as the PDE and SDE cases. Through these examples, we show the idea of EnSIL works quite well even under complicated situations involving both temporal and spatial changes. 

\subsection{Michaelis–Menten kinetics}
In 1913, to explain how enzymes accelerate a chemical reaction, Michaelis and Menten proposed a simple model by including both enzyme binding to and unbinding from a substrate \cite{michaelis1913kinetik}. Its reversible version reads 
\begin{equation}\label{Michaelis_Menten}
E+S 
\underset{k_{1}^-}{\stackrel{k_{1}^+}{\rightleftarrows}} E S \underset{k_{2}^-}{\stackrel{k_{2}^+}{\rightleftarrows}} E+P,
\end{equation}
where $E$ stands for the enzyme, $S$ for the substrate, $E S$ for a complex of enzyme and substrate, and $P$ for the product. Clearly, after two sequential forward reactions, a substrate is transformed into a product, meanwhile the enzyme keeps unchanged. $k$'s are reaction rate constants. 

According to the law of mass action, a number of ODEs for reactants and products are arrived at,
\begin{subequations}
\begin{align}
&\frac{d[E]}{d t}=-k_{1}^+[E][S]+k_{1}^-[E S]+k_{2}^+[E S]-k_{2}^-[E][P], \\
&\frac{d[S]}{d t}=-k_{1}^+[E][S]+k_{1}^-[E S],\label{MMa} \\
&\frac{d[E S]}{d t}=k_{1}^+[E][S]-k_{1}^-[E S]-k_{2}^+[E S]+k_{2}^-[E][P], \\
&\frac{d[P]}{d t}=k_{2}^+[E S]-k_{2}^-[E][P],\label{MMb}
\end{align}
\label{MMeq}
\end{subequations}
where the rate constants are set as $k_{1}^+=1, k_{1}^-=10, k_{2}^+=1, k_{2}^-=0$. And the initial conditions are $[E](0)=20$, $[S](0)=50$, $[ES](0)=10$ and $[P](0)=10$ respectively.

\begin{figure}[h]
\centering
\includegraphics[width=1.0\linewidth]{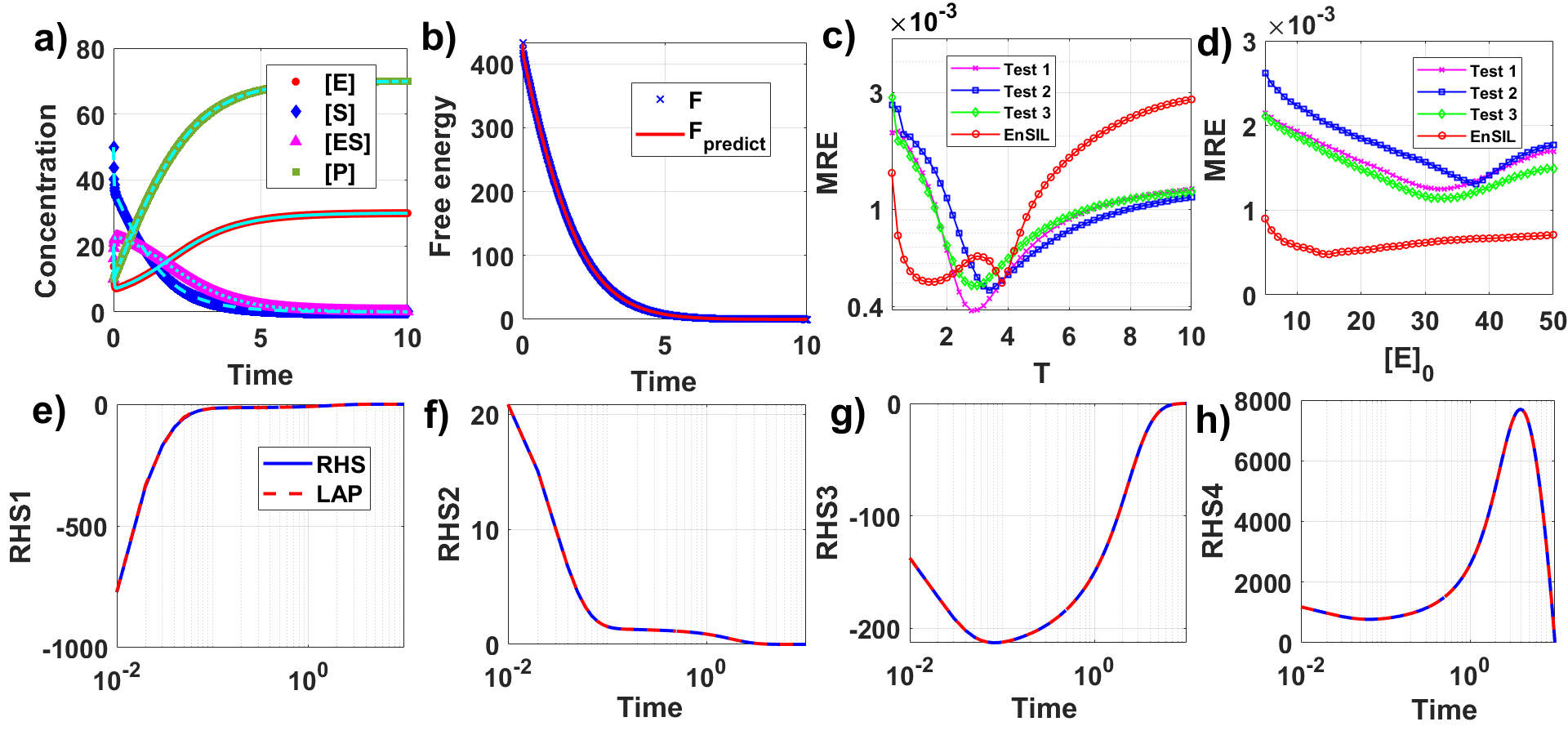}
\caption{Predictions of the EnSIL method on MM kinetics. (a-b) Comparison of predictions of the EnSIL method and exact solutions on species concentration and free energy. (c-d) Mean relative errors of the EnSIL method and other three tests with respect to a varied reference time $T$ and initial concentration of enzyme $[E]_0$. (e-h). Accuracy of the Legendre polynomial approximation on the right-hand-side terms in Eq. \ref{entropy_MM}.}
\label{fig.MM}
\end{figure}

On the basis of Eqs. \eqref{REF} and \eqref{relative_energy},
the free energy dissipation rate $d F(t) / d t$ is given by
\begin{equation}
\begin{aligned}
\frac{d F(t)}{d t}  & =k_1^+[E][S]\ln{A}-k_{1}^-[ES]\ln{A}+k_{2}^+[ES]\ln{B}-k_{2}^-[E][P]\ln{B}
=\mathbf{\mathcal{P}}\cdot\mathbf{\Theta},
\end{aligned}
\label{entropy_MM}
\end{equation}
where $A = \displaystyle\frac{[E]^{s s}[S]^{s s}[E S]}{[E][S][E S]^{s s}}, B = \displaystyle\frac{[E][P][E S]^{s s}}{[E]^{s s}[P]^{s s}[E S]}$. $\mathbf{\mathcal{P}}=(P_1, P_2, P_3, P4)$ represent Legendre polynomials and $\mathbf{\Theta} = (k_1,k_{-1},k_2,k_{-2})^T$ is the coefficient matrix. Referring to above formula, we can implement the EnSIL algorithm to study the inverse problem of MM reaction. Under a setup of the time step $\Delta t = 0.001$ and the simulation time  $t\in[0,10]$, our method shows higher accuracy in discovering the unknown rate constants (mean relative error is 0.41\%) than SINDy (see Table. S1). The corresponding predicted trajectories by EnSIL are compared with the true ones in Figs. \ref{fig.MM}(a-b). In addition, Figs.\ref{fig.MM}(e-h) provide a detailed examination on the effectiveness of Legendre polynomial approximation for terms on the right-hand side of Eq. \ref{entropy_MM}. 

To compare the performance of EnSIL and other classical machine learning algorithms, here we design four different tests. (i) Test 1: Integral SINDy algorithm for  Eq.\ref{MMa} and Eq.\ref{MMb}; (ii) Test 2: a direct regression of all four equations in Eq. \ref{MMeq} by Integral SINDy ; (iii) Test 3: Integral SINDy  algorithm for Eq. \ref{MMeq} plus Eq. \ref{entropy_MM}; In this case, the entropy balance equation has been used as an additional constraint. (4) Test 4: Entropy Informed Learning based on Eq. \ref{entropy_MM} solely. A few remarks follow Test 1. According to the law of mass conservation, we have $[E] + [ES] = const, [S] + [P] + [ES]= const$. It means there are only two independent variables, saying $[S]$ and $[P]$, which we can utilize to accurately deduce the original whole model. 

In principle, a calculation of the Helmholtz free energy requires knowledge of the concentrations at the equilibrium state, which is quite challenging since it relates to the infinite time. However, in practice we can safely replace the equilibrium concentration with those at large enough reference time point $T$, and still keep the basic mathematical properties of the new free energy unchanged till time $T$. This is a direct consequence of the Kullback-Leibler divergence and provides a lot of facilities in real applications of EnSIL.

In Fig.\ref{fig.MM}(c-d), the influence of varied reference time $T$ and initial concentrations of substance $[E](0)$ have been compared among four different tests. When the reference time $T\leq 2.2$, the EnSIL method outperforms all others. The highest accuracy is achieved by Test 2 when $T\geq 3.9$. Furthermore, the EnSIL method has the lowest mean relative errors in deducing model parameters, regardless of how we alter the initial concentrations.

\subsection{Schlögl model}

Schlögl model is a remarkable chemical reaction network which  exhibits bistability  \cite{schlogl1972chemical,vellela2009stochastic}. This reaction depicts a conversion between $A$ and $B$ through an intermediate species $X$,
\begin{equation}\label{Schlogl}
\begin{gathered}
2X + A
\underset{k_{2}}{\stackrel{k_{1}}{\rightleftarrows}} 3X,\\
B
\underset{k_{4}}{\stackrel{k_{3}}{\rightleftarrows}} X.
\end{gathered}
\end{equation}
where $X,A,B$ are three reactants. $x=[X]$, $a = [A]$, $b=[B]$ denote the concentration of each reactant, while $k_i\; (i = 1,2,3,4)$ are rate constants. The deterministic governing equations are given through
\begin{equation}
\begin{aligned}
&\frac{d x}{d t}=k_{1} x^{2} a-k_{2} x^{3}+k_{3}b-k_{4}x, \\
&\frac{d a}{d t}=-k_{1} x^{2} a+k_{2} x^{3}, \\
&\frac{d b}{d t}=-k_{3} b+k_{4} x.
\end{aligned}
\label{Schlogl-mae}
\end{equation}
The free energy dissipation rate $d F(t)/dt$ becomes 
\begin{equation}
\label{Schlogl-ent}
       \frac{d F}{d t}=\left(k_2 x^3-k_1 x^2 a\right) \ln \frac{a x^{s s}}{x a^{s s}}+\left(k_4 x-k_3 b\right) \ln \frac{b x^{s s}}{x b^{s s}}.
\end{equation}

\begin{figure}[h]
\centering
\includegraphics[width=1.0\linewidth]{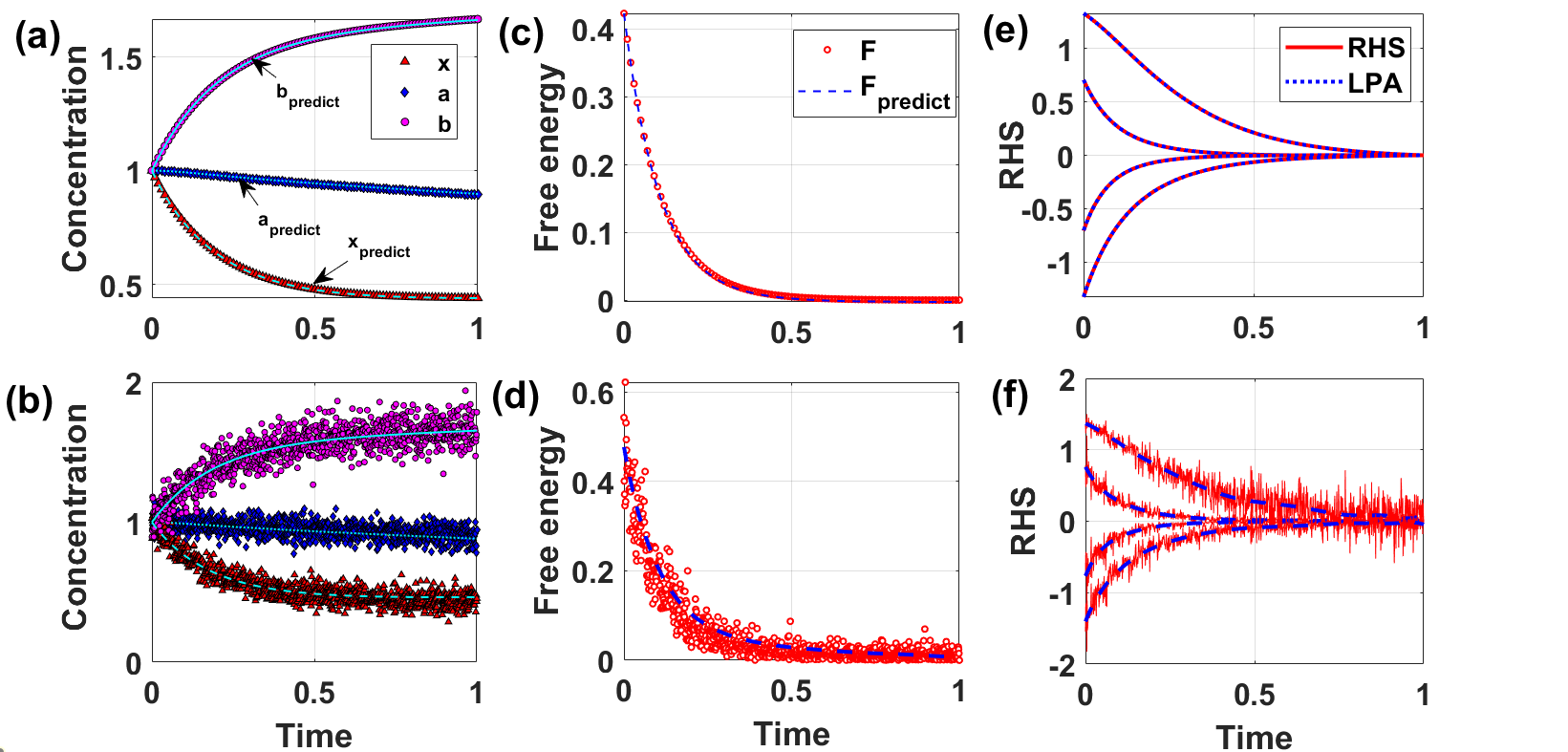}
\caption{Performance of the EnSIL method on the Schlögl model in the presence of (a,c) noiseless data and (b,d) noisy data (combined with $5\%$ white noise). Legendre polynomials are adopted to fit four right-hand-side terms in Eq. \ref{Schlogl-ent} with (e) noiseless and (f) noisy data, respectively.}
\label{fig.Schlogl}
\end{figure}

In this case, the EnSIL method is adopted for both noiseless and noisy data (see Figs. \ref{fig.Schlogl}(a-b)). Additive white noise $\eta_j(t)=\epsilon\|\boldsymbol c_j\|_{\infty} \boldsymbol{\xi},\; \boldsymbol{\xi} \sim \mathcal{N}\left(0, \sigma^2\right)\; (j=1,2,3)$ 
with noise intensity $\epsilon$ ranging from $0.001\%$ to $10\%$ mixed with the noiseless concentration data have been examined. From Table \ref{Tab.schlogl}, it can be seen that EnSIL is more reliable to the noisy data than SINDy. Time steps also play an important role in the  accuracy of machine learning algorithms \cite{HU2022111203}. Generally speaking, the integration-based approach is significantly more tolerant than differentiation-based approaches for large time steps (see more details in Table \ref{Tab.schlogl}).  Further compared to directly learning the mass-action equations in Eq. \ref{Schlogl-mae}, EnSIL based on learning the entropy balance equation in Eq. \ref{Schlogl-ent} is more accurate under various time steps.


\begin{table}[h]
\centering
\setlength\extrarowheight{10pt}
\setlength{\tabcolsep}{6pt}
\begin{tabular}{|ccc|cccc|c|}
\hline
\multicolumn{3}{|c|}{\multirow{2}{*}{Schlögl model}}                                                                                                        & \multicolumn{4}{c|}{Parameters}                                                                  & \multirow{2}{*}{MRE} \\ \cline{4-7}
\multicolumn{3}{|c|}{}                                                                                                                                      & \multicolumn{1}{c|}{$k_1$}  & \multicolumn{1}{c|}{$k_2$}  & \multicolumn{1}{c|}{$k_3$}  & $k_4$  &                      \\ \hline
\multicolumn{3}{|c|}{True}                                                                                                                                  & \multicolumn{1}{c|}{1}      & \multicolumn{1}{c|}{1}      & \multicolumn{1}{c|}{1}      & 4      & -                    \\ \hline
\multicolumn{1}{|c|}{\multirow{8}{*}{\begin{tabular}[c]{@{}c@{}}Time step\\ $\Delta t$\end{tabular}}} & \multicolumn{1}{c|}{EnSIL}    & $\Delta t=0.001$ & \multicolumn{1}{c|}{0.9972} & \multicolumn{1}{c|}{0.9968} & \multicolumn{1}{c|}{1.0004} & 4.0085 & 0.18\%               \\ \cline{2-8} 
\multicolumn{1}{|c|}{}                                                                                & \multicolumn{1}{c|}{SINDy-based} & $\Delta t=0.001$ & \multicolumn{1}{c|}{0.9831} & \multicolumn{1}{c|}{0.9768} & \multicolumn{1}{c|}{0.9960} & 3.9866 & 1.19\%               \\ \cline{2-8} 
\multicolumn{1}{|c|}{}                                                                                & \multicolumn{1}{c|}{EnSIL}    & $\Delta t=0.01$  & \multicolumn{1}{c|}{0.9977} & \multicolumn{1}{c|}{0.9976} & \multicolumn{1}{c|}{1.0004} & 4.0080 & 0.15\%               \\ \cline{2-8} 
\multicolumn{1}{|c|}{}                                                                                & \multicolumn{1}{c|}{SINDy-based} & $\Delta t=0.01$  & \multicolumn{1}{c|}{0.9307} & \multicolumn{1}{c|}{0.9063} & \multicolumn{1}{c|}{0.9777} & 3.9135 & 5.17\%               \\ \cline{2-8} 
\multicolumn{1}{|c|}{}                                                                                & \multicolumn{1}{c|}{EnSIL}    & $\Delta t=0.1$   & \multicolumn{1}{c|}{0.9943} & \multicolumn{1}{c|}{0.9929} & \multicolumn{1}{c|}{0.9975} & 4.0011 & 0.31\%               \\ \cline{2-8} 
\multicolumn{1}{|c|}{}                                                                                & \multicolumn{1}{c|}{SINDy-based} & $\Delta t=0.1$   & \multicolumn{1}{c|}{0.5751} & \multicolumn{1}{c|}{0.5459} & \multicolumn{1}{c|}{0.7858} & 3.1476 & 32.66\%              \\ \cline{2-8} 
\multicolumn{1}{|c|}{}                                                                                & \multicolumn{1}{c|}{EnSIL}    & $\Delta t=0.5$   & \multicolumn{1}{c|}{1.0019} & \multicolumn{1}{c|}{1.0094} & \multicolumn{1}{c|}{0.8958} & 3.6181 & 4.22\%               \\ \cline{2-8} 
\multicolumn{1}{|c|}{}                                                                                & \multicolumn{1}{c|}{SINDy-based} & $\Delta t=0.5$   & \multicolumn{1}{c|}{0.2343} & \multicolumn{1}{c|}{0.2246} & \multicolumn{1}{c|}{0.3693} & 1.4793 & 70.01\%              \\ \hline
\multicolumn{1}{|c|}{\multirow{6}{*}{Noise intensity}}                                                & \multicolumn{1}{c|}{EnSIL}    & 0.1\%            & \multicolumn{1}{c|}{0.9996} & \multicolumn{1}{c|}{0.9978} & \multicolumn{1}{c|}{1.0050} & 4.0142 & 0.24\%               \\ \cline{2-8} 
\multicolumn{1}{|c|}{}                                                                                & \multicolumn{1}{c|}{SINDy-based} & 0.1\%            & \multicolumn{1}{c|}{0.9540} & \multicolumn{1}{c|}{0.9426} & \multicolumn{1}{c|}{0.9928} & 3.9706 & 2.95\%               \\ \cline{2-8} 
\multicolumn{1}{|c|}{}                                                                                & \multicolumn{1}{c|}{EnSIL}    & 1\%              & \multicolumn{1}{c|}{0.9882} & \multicolumn{1}{c|}{0.9825} & \multicolumn{1}{c|}{1.0178} & 4.0193 & 1.12\%               \\ \cline{2-8} 
\multicolumn{1}{|c|}{}                                                                                & \multicolumn{1}{c|}{SINDy-based} & 1\%              & \multicolumn{1}{c|}{1.0087} & \multicolumn{1}{c|}{1.0615} & \multicolumn{1}{c|}{0.9022} & 3.7234 & 5.93\%               \\ \cline{2-8} 
\multicolumn{1}{|c|}{}                                                                                & \multicolumn{1}{c|}{EnSIL}    & 5\%              & \multicolumn{1}{c|}{0.9318} & \multicolumn{1}{c|}{0.9072} & \multicolumn{1}{c|}{1.0744} & 4.0405 & 5.22\%               \\ \cline{2-8} 
\multicolumn{1}{|c|}{}                                                                                & \multicolumn{1}{c|}{SINDy-based} & 5\%              & \multicolumn{1}{c|}{2.2134} & \multicolumn{1}{c|}{1.5515} & \multicolumn{1}{c|}{1.0658} & 4.3779 & -                    \\ \hline
\end{tabular}
\caption{Impact of the time step $\Delta t$ and noise intensity on the precision of EnSIL and SINDy methods for Schlogl model. In the time-step test no noise is added, while $\Delta t$ is fixed at 0.01 for the noise-intensity test.}
\label{Tab.schlogl}
\end{table}


\subsection{Chemical Lorenz equation}
The MM reaction and Schlogl model are both closed systems since they are free from species changing with the reservoirs.  In this example, We proceed to consider an open nonlinear chemical reaction system, in which concentrations of some certain species, called chemostat, are restrained to external control. 
This nonlinear reaction system is chaotic since its evolution is highly sensitive to the initial conditions. And it brings significant challenges to the machine-learning-based discoveries of differential equations from the time-series data since a small deviation in the initial data may produce totally different trajectories.  The chemical Lorenz equations \cite{zhang2020dynamic} are one such system which offers a critical test on our EnSIL method. It was established to imitate the classic Lorenz equation by constructing an artificial reaction scheme,
\begin{equation}\label{MM}
\begin{aligned}
& R_1+X_1+X_2 \stackrel{k_1}{\longrightarrow} 2 X_1+X_2, \quad
R_2+X_1+X_2 \stackrel{k_2}{\longrightarrow} X_1+2 X_2, \quad R_3+X_3 \stackrel{k_3}{\longrightarrow} 2 X_3, \\
& X_1+X_2+X_3 \stackrel{k_4}{\longrightarrow} X_1+2 X_3, \quad X_2+X_3 \stackrel{k_5}{\longrightarrow} 2 X_2, \quad 2 X_1 \stackrel{k_6}{\longrightarrow} P_1, \\
& 2 X_2 \stackrel{k_7}{\longrightarrow} P_2, \quad X_2 \stackrel{k_8}{\longrightarrow} P_3,  \quad X_1+X_3 \stackrel{k_9}{\longrightarrow} X_1+P_4, \quad 2 X_3 \stackrel{k_{10}}{\longrightarrow} P_5 .
\end{aligned}
\end{equation}
where $R_1,R_2,R_3$ are three reactant reservoirs, whose concentrations are remained constant. Based on the law of mass-action, the governing equations \cite{zhang2020dynamic} read 
\begin{equation}
\begin{aligned}
&\frac{d x_{1}}{d t}=k_{1} r_{1} x_{1} x_{2}-2 k_{6} x_{1}^{2}, \\
&\frac{d x_{2}}{d t}=k_{2} r_{2} x_{1} x_{2}-k_{4} x_{1} x_{2} x_{3}+k_{5} x_{2} x_{3}-2 k_{7} x_{2}^{2}-k_{8} x_{2}, \\
&\frac{d x_{3}}{d t}=k_{3} r_{3} x_{3}+k_{4} x_{1} x_{2} x_{3}-k_{5} x_{2} x_{3}-k_{9} x_{1} x_{3}-2 k_{10} x_{3}^{2},
\end{aligned}
\end{equation}
where the reaction constants are set $k_{1}=0.001, k_{2}=1, k_{3}=10000, k_{4}=0.0001, k_{5}=1$,
$k_{6}=0.05, k_{7}=0.005, k_{8}=9900, k_{9}=1, k_{10}=0.0133$. $[R_i]=r_i$ denotes the concentration for the reactant reservoir $R_i~(i=1,2,3)$, and $[X_{j}]=x_{j},~(j=1,2,3)$ for concentrations of reactants. 

As to the chemical Lorenz equations, since it is an open system, the entropy function $S=\sum_{i=1}^3 x_i\ln (x_i)-x_i$ is adopted instead of previous free energy function, whose change rate is given by
\begin{equation}\label{SLorenz}
\begin{aligned}
\frac{d S}{d t} =&k_{1} r_{1} x_{1} x_{2} \ln x_1+k_{2} r_{2} x_{1} x_{2} \ln x_2 + k_{3} r_{3} x_{3} \ln x_3+k_{4} x_{1} x_{2} x_{3}\ln (x_3/x_2) +k_{5} x_{2} x_{3}\ln (x_2/x_3)\\
&-2 k_{6} x_{1}^{2} \ln x_1 -2 k_{7} x_{2}^{2} \ln x_2 - k_{8} x_{2} \ln x_2 - k_{9} x_{1} x_{3} \ln x_3 -2 k_{10} x_{3}^{2} \ln x_3.
\end{aligned}
\end{equation}

\begin{figure}[h]
\centering  
\includegraphics[width=1.0\linewidth]{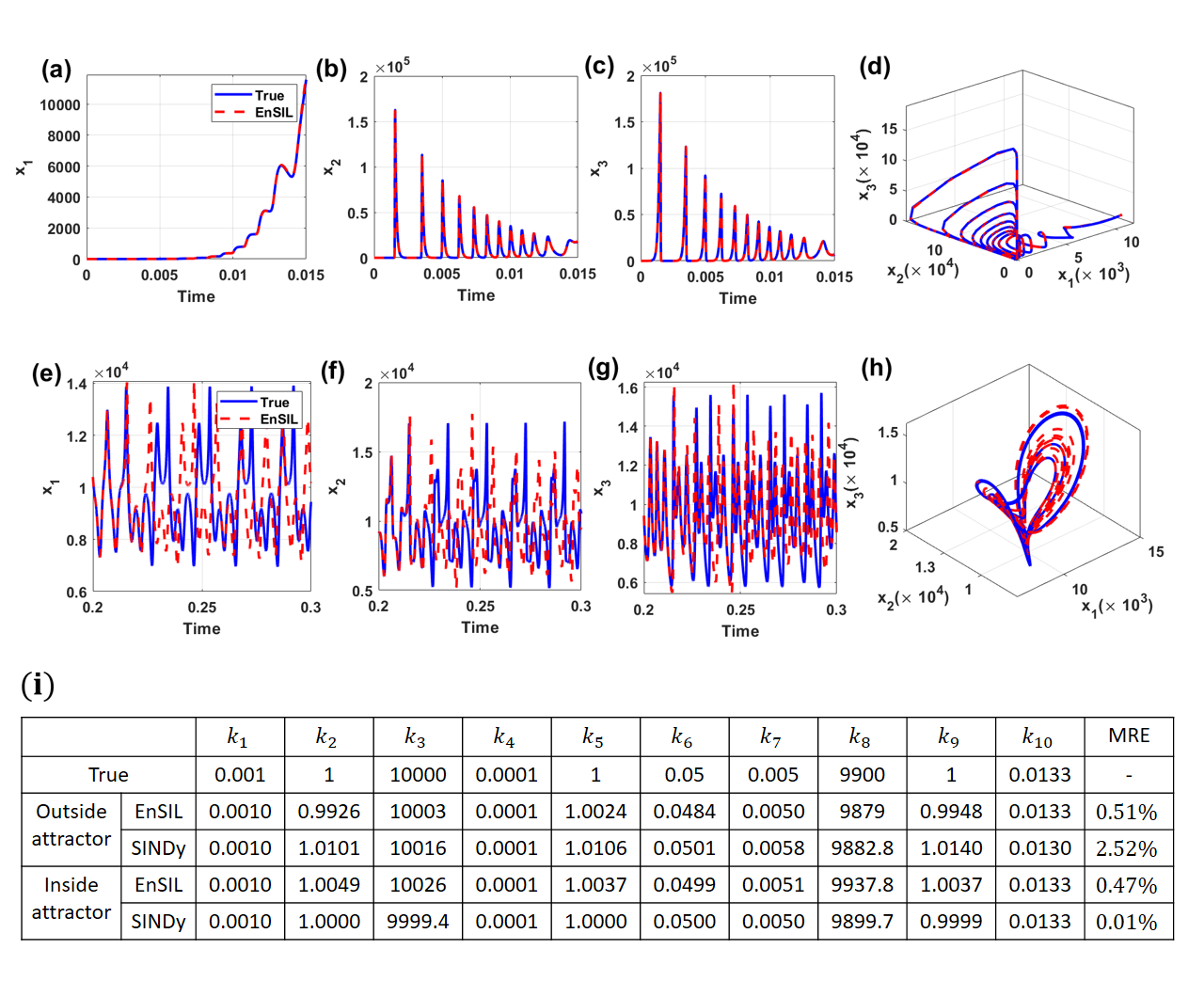}
\caption{Predicted time trajectories of $x_1$, $x_2$, $x_3$ in the chemical Lorenz equation by the EnSIL method (a-d) outside and (e-h) inside the strange attractor. Blue solid lines stand for numerical data, while red dashed lines for predictions of EnSIL method. (i) The predicted values of ten coefficients in the chemical Lorenz equation are compared for two methods.}
\label{fig.lorenz}
\end{figure}

Here we explore the inverse problem under two distinct situations: before entering the Lorenz attractor and evolution within the Lorentz attractor. In the first case (Figs. \ref{fig.lorenz}(a-d)), the EnSIL shows a nearly perfect performance and correctly derives all ten unknown coefficients that vary over eight orders of magnitude. While in the second case (Figs. \ref{fig.lorenz}(e-h)), since the dynamics is highly chaotic, the entropy evolution is chaotic too. The mean relative errors of EnSIL are a bit higher (about $0.47\%$) than SINDy, but are still within a quite satisfactory region.
                     

\subsection{Porous medium equation}
After the studies on ODEs, now we turn to far more complicated PDEs. The porous medium equation (PME), also known as a nonlinear parabolic equation, takes the following form,
\begin{equation}
\frac{\partial \rho}{\partial t}=\sum_{i=1}^d\sum_{j=1}^d D_{ij}\frac{\partial^2}{\partial x_i\partial x_j}(\rho^m),
\end{equation}
where $\rho=\rho(\mathbf{x},t): \mathbb{R}^d \rightarrow \mathbb{R}$ is a time–dependent density function and $(D_{ij})\geq0$ is a semi-positive definite matrix. The exponent satisfies $m > 1 - \displaystyle\frac{1}{d}$ and $m> \displaystyle\frac{d}{d+2}$ according to Ref. \cite{otto2001geometry}. And an entropy functional for PME is given by,
\begin{equation}
S(\rho)=\left\{\begin{array}{l}
\displaystyle\frac{1}{m-1} \int \rho^m, m \neq 1, \\
\\
\displaystyle\int \rho \ln \rho, m=1.
\end{array}\right.
\end{equation}

\begin{figure}[h]
\centering  
\includegraphics[width=1.0\linewidth]{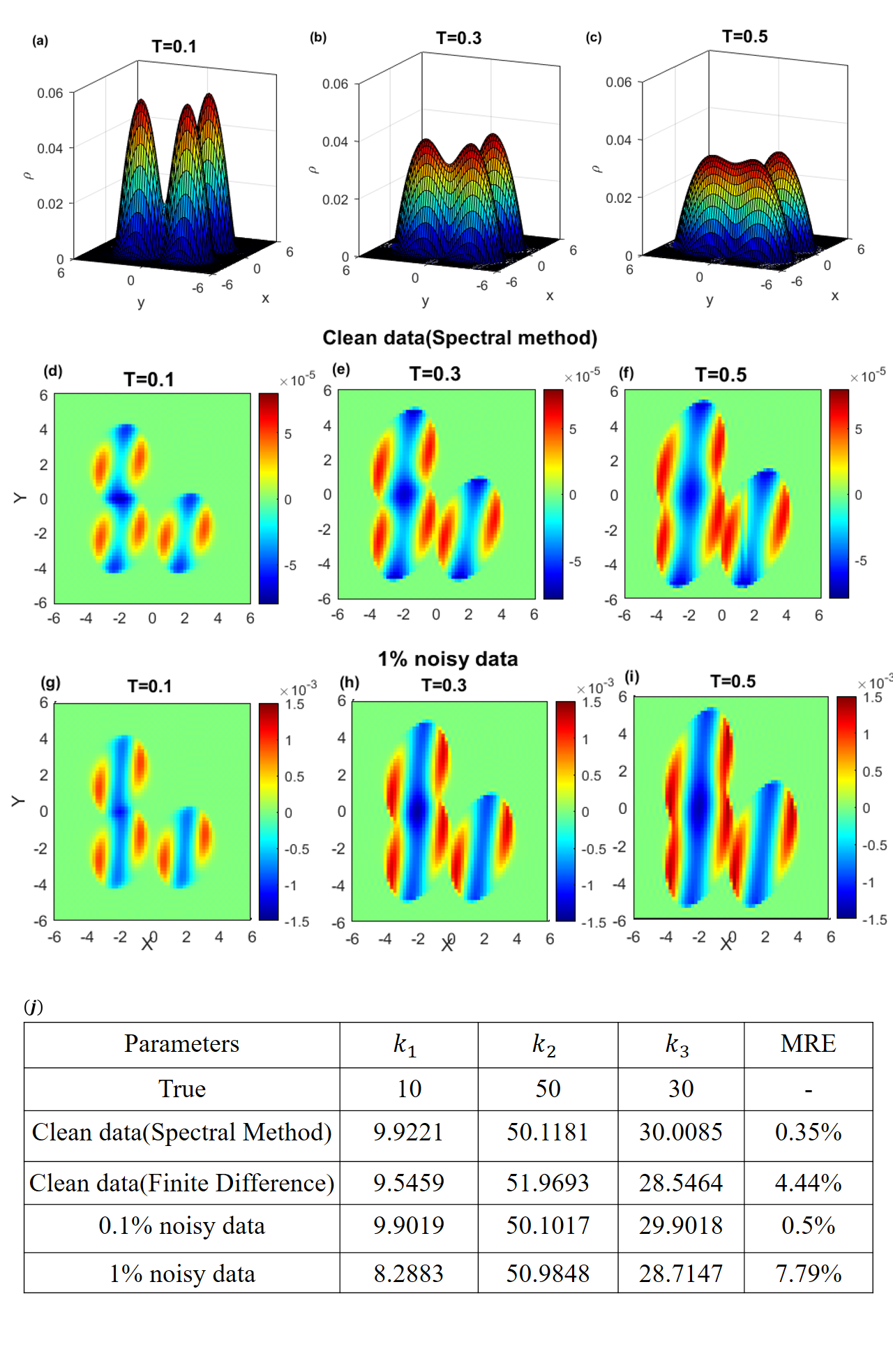}
\caption{Solutions to the porous medium equation at time (a) $t=0$ (b) $t=0.3$ and (c) $t=0.5$. The initial distribution is set as $\rho(x,y,0)=\frac{1}{3\pi}(e^{-(x+2)^2-(y+2)^2}+e^{-(x+2)^2-(y-2)^2}+e^{-(x-2)^2-(y+2)^2})$. The relative errors of predictions of the EnSIL method are highlighted in (d-f) for noiseless data and (g-i) for noisy data. The spectral method is adopted to improve the accuracy. (j) Summary on the accuracy of the EnSIL method for solving the inverse problem of the 2D porous medium equation.}
\label{fig.pme}
\end{figure}

Here we test our EnSIL method on a two-dimensional case, i.e.
\begin{equation}\label{pme2}
\frac{\partial \rho}{\partial t}-k_1(\rho^m)_{xx}-k_2(\rho^m)_{yy}-k_3(\rho^m)_{xy}=0
\end{equation}
where $m = 2$ and $k_1 = 10, k_2 = 50, k_3 = 30$. The initial distribution of this case is shown in Fig. \ref{fig.pme}(a). The time step size is set as $\Delta t=0.01$ and the total simulation time is from 0 to 0.5. And the spatial domain is discretized into 128$\times 128$ grids. Besides, white noise is added to the clean data by multiplying $1+\epsilon$, where $\epsilon$ obeys a standard normal distribution. In the light of results presented in Section 2.3, the entropy change rate ${d S}/{d t}$ reads,
\begin{equation}
\frac{d S}{d t} =-m^2\left(k_1 \int \rho^{2 m-3}\left(\rho_x\right)^2+k_2 \int \rho^{2 m-3}\left(\rho_y\right)^2+k_3 \int \rho^{2 m-3} \rho_x \rho_y\right).
\end{equation}

For PME, the spectral method is adopted to calculate the spatial derivatives, i.e.
\begin{equation*}
0=\frac{\partial \rho}{\partial t}-\mathcal{F}^{-1}\left(\mathcal{F}\left(\sum_{i=1}^d\sum_{j=1}^d D_{ij}\frac{\partial^2}{\partial x_i\partial x_j}(\rho^m)\right)\right)=\frac{\partial \rho}{\partial t}+\mathcal{F}^{-1}\left(\sum_{i=1}^d\sum_{j=1}^d \omega_i\omega_jD_{ij}\mathcal{F}\left(\rho^m\right)\right).
\end{equation*}
in which a function's derivative is represented via the Fourier transform $\mathcal{F}$,
\begin{equation*}
\mathcal{F}\left(f^{\prime}(x)\right) =\int_{-\infty}^{\infty} f^{\prime}(x) e^{-i \omega x} d x
=\left[f^{\prime}(x) e^{-i \omega x}\right]_{-\infty}^{\infty}-\int_{-\infty}^{\infty} f(x)\left[-i \omega e^{-i \omega x}\right] d x=i \omega \mathcal{F}(f(x)).
\end{equation*}
Fig. \ref{fig.pme} summarizes the ground truth and  learned parameters by using the EnSIL method in the presence of external noises or not. The spectral method shows higher accuracy than the finite-difference schemes.

\subsection{Stochastic dynamics with a double-well potential}
\begin{figure}
\centering
\includegraphics[width=1.0\linewidth]{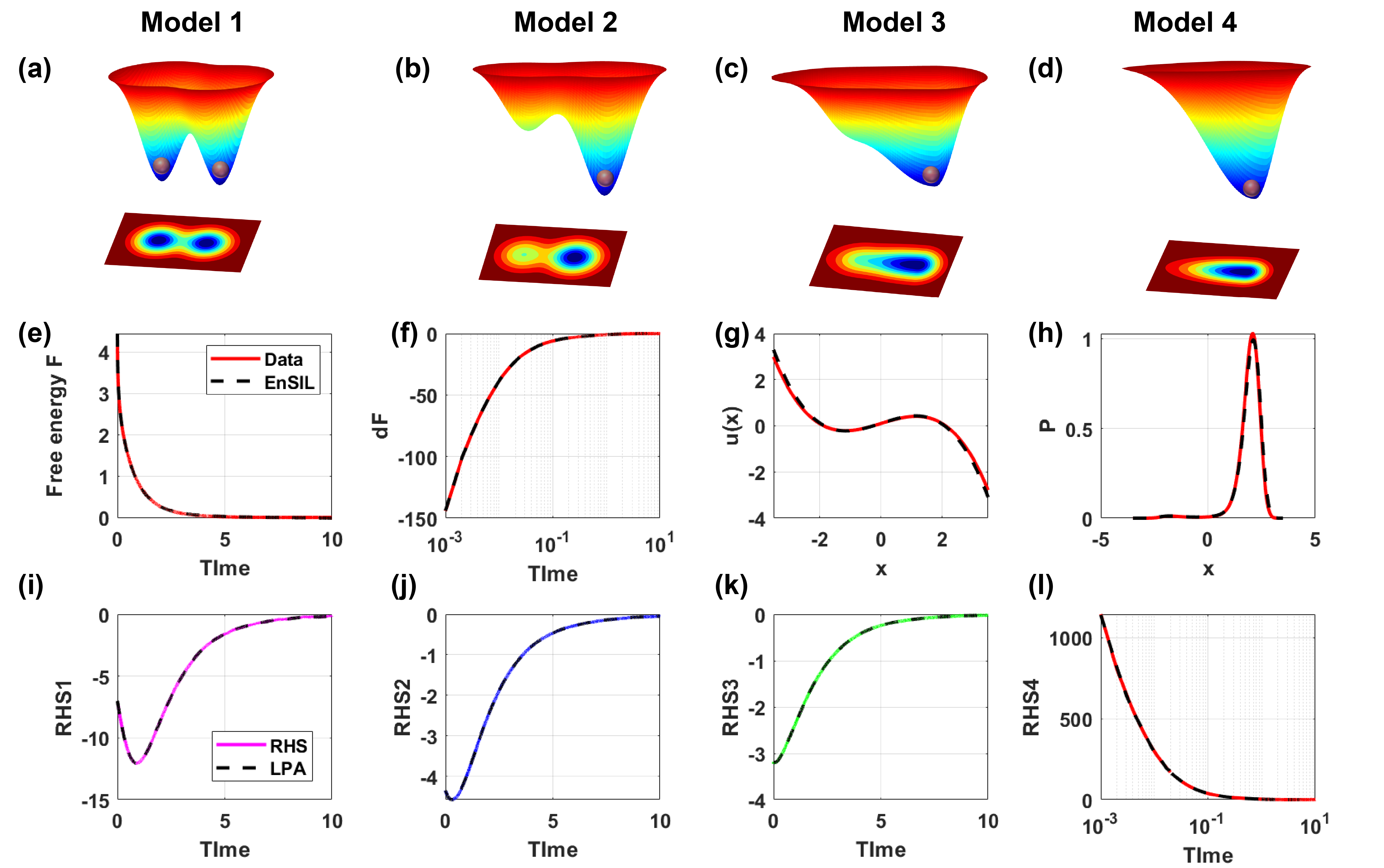}
\caption{(a-d) Schematic diagram of the potential energy surface used for the stochastic double-well systems. As an illustration, the predicted (e) free energy, (f) free energy dissipation rate, (g) drift term and (h) probability density for Model 2 at time T = 10 are compared with their exact values, respectively. The accuracy of Legendre polynomial approximation is highlighted in (i-l).}
\label{fig.doublewell}
\end{figure}

Bistability, bifurcation, chaos and hysteresis are typical phenomena of nonlinear systems. Consider a stochastic Langevin-type dynamics \cite{nolte2015introduction}
\begin{equation}
d x=-\nabla_x V(x) d t+\sigma d B_t=u(x)dt+\sigma dB_t,\quad x \in \mathbb{R}, 
\label{SDE}
\end{equation}
with a double-well potential 
\begin{equation}\label{potential}
	V(x)=\frac{a}{4} x^{4}-\frac{b}{2} x^{2}-c x,
\end{equation}
where $c$ is the bias parameter for the system. $B_{t}$ denotes the standard Brownian motion. Fix $\sigma=0.5$ and $a=0.1$. The numerical setups in Table \ref{Tab.SDE} and Fig. \ref{fig.doublewell}(a-d) involve both symmetric ($c = 0$) and asymmetric ($c \neq 0$) double-wells.

The bias parameter $c$ has a significant impact not only on the energy difference between the left and right potential wells, but also on the dynamical behaviors of the system. When $c = 0$ the double well has two identical energy minimums, While $c>0$ a single global energy minimum is found in the right-hand well. In particular, the potential function $V(x)$ has three fixed points when $|c|<c^{*}=0.3079$. The saddle point is in between the two stable fixed points, which correspond to the double-well potential surface of $V(x)$ (see Fig.S1(d)).

Based on the results presented in Sec. 2.4, the inverse problem of SDE is converted into an equivalent inverse problem of PDE in a sense of Itô's calculus -- learning the unknown coefficients of a Fokker-Planck equation. To do this, we sample $10^6$ stochastic trajectories through solving the SDE in Eq. \ref{SDE} by the Euler–Maruyama method. And a non-parametric estimation method -- the Kernel Density Estimation (KDE) is adopted to obtain a smooth probability density function. For low-dimensional problems, KDE yields very accurate estimations of the probability density function. To approximate the probability density function in a high-dimensional space, Chen et al. \cite{chen2021solving} have discussed some advanced algorithms. Ma et al. \cite{ma2021learning} applied the weak form of a Fokker-Planck equation to avoid a direct calculation of the probability density function.

As summarized in Table \ref{Tab.SDE}, we compare EnSIL with the PDE-FIND algorithm in solving the inverse problem of SDE. And it is found that in all four cases the EnSIL method is much more efficient and accurate in finding out the correct different combinations of drift and diffusion terms than PDE-FIND.

\begin{table}[H]
\setlength\extrarowheight{12pt}
\setlength{\tabcolsep}{10pt}
\begin{tabular}{|cc|ccc|c|c|}
\hline
\multicolumn{2}{|l|}{\multirow{2}{*}{Double-well potential stochastic system}}             & \multicolumn{3}{c|}{$u(x)$}                                                              & $D$                         & \multicolumn{1}{l|}{\multirow{2}{*}{MRE}} \\ \cline{3-6}
\multicolumn{2}{|l|}{}                                            & \multicolumn{1}{c|}{$x^3$}   & \multicolumn{1}{c|}{$x^1$}  & 1                           & 1                           & \multicolumn{1}{l|}{}                     \\ \hline
\multicolumn{1}{|l|}{\multirow{3}{*}{Model 1}} & True             & \multicolumn{1}{c|}{-0.1000} & \multicolumn{1}{c|}{0.3000} & 0                           & 0.5000                      & -                                         \\ \cline{2-7} 
\multicolumn{1}{|l|}{}                         & EnSIL & \multicolumn{1}{c|}{ -0.1033} & \multicolumn{1}{c|}{0.3059} & 0                           & 0.5050                      & 2.09\%                                    \\ \cline{2-7} 
\multicolumn{1}{|l|}{}                         & PDE-FIND            & \multicolumn{1}{c|}{-0.0746} & \multicolumn{1}{c|}{ 0.2647} & 0                           & 0.4271                      & 14.17\%                                   \\ \hline
\multicolumn{1}{|c|}{\multirow{3}{*}{Model 2}} & True             & \multicolumn{1}{c|}{-0.1000} & \multicolumn{1}{c|}{0.4000} & 0.1                         & 0.5000                      & -                                         \\ \cline{2-7} 
\multicolumn{1}{|c|}{}                         & EnSIL & \multicolumn{1}{c|}{-0.1093} & \multicolumn{1}{c|}{ 0.4243} &  0.1041                     & 0.4994                     & 4.93\%                                    \\ \cline{2-7} 
\multicolumn{1}{|c|}{}                         & PDE-FIND            & \multicolumn{1}{c|}{-0.0391} & \multicolumn{1}{c|}{0.1832} &  0.0333                      &  0.2637                      & 50.61\%                                   \\ \hline
\multicolumn{1}{|c|}{\multirow{3}{*}{Model 3}} & True             & \multicolumn{1}{c|}{-0.1000} & \multicolumn{1}{c|}{0.4000} & 0.3                         & 0.5000                      & -                                         \\ \cline{2-7} 
\multicolumn{1}{|c|}{}                         & EnSIL & \multicolumn{1}{c|}{-0.1060} & \multicolumn{1}{c|}{0.4239} & 0.2892                     & 0.5261                      & 5.20\%                                    \\ \cline{2-7} 
\multicolumn{1}{|c|}{}                         & PDE-FIND            & \multicolumn{1}{l|}{-0.1122} & \multicolumn{1}{l|}{0.1335} & \multicolumn{1}{l|}{0.4689} & \multicolumn{1}{l|}{0.4136} & \multicolumn{1}{c|}{39.55\%}              \\ \hline
\multicolumn{1}{|c|}{\multirow{3}{*}{Model 4}} & True             & \multicolumn{1}{c|}{-0.1000} & \multicolumn{1}{c|}{0.4000} & 0.35                        & 0.5000                      & -                                         \\ \cline{2-7} 
\multicolumn{1}{|c|}{}                         & EnSIL & \multicolumn{1}{l|}{-0.1065} & \multicolumn{1}{l|}{0.3976} & \multicolumn{1}{l|}{0.3763} & \multicolumn{1}{l|}{0.5115} & \multicolumn{1}{c|}{4.59\%}               \\ \cline{2-7} 
\multicolumn{1}{|c|}{}                         & PDE-FIND            & \multicolumn{1}{l|}{-0.0828} & \multicolumn{1}{l|}{0.3451} & \multicolumn{1}{l|}{0.3080} & \multicolumn{1}{l|}{0.4161} & \multicolumn{1}{c|}{14.81\%}              \\ \hline
\end{tabular}
\caption{Performance of the EnSIL method for solving the inverse problem of the double-well potential stochastic system. The number of samples is $10^6$.}
\label{Tab.SDE}
\end{table}

\section{Conclusion and discussion}
In this paper, we have proposed a novel entropy-structure informed method to solve the inverse problems of XDEs. Motivated by the elegant entropy structure for various equations, like the existence of a Lyapunov function for a number of ODEs, the asymptotic stability of solutions to the time evolutionary equations characterized by the entropy or free energy function, the entropy condition for hyperbolic equations in the presence of shocks, etc., we explore the possibility of solving the inverse problems of XDEs  (including ODEs, PDEs, and SDEs) in an alternative way. By transforming the original XDEs into a single entropy balance equation, we show that adequate information on the unknown coefficients of XDEs could still be extracted for a number of examples under study. Further combining with Legendre polynomial approximation and integral-based regression methods, our entropy-structure informed method could achieve comparable performance with other classical machine learning algorithms, like SINDy and PDE-FIND.

The accuracy and robustness of our EnSIL method are assessed on a number of classical examples, including chemical reaction kinetics in both closed and open systems, the 2D porous medium equation, and the Langevin-type stochastic dynamics with four different kinds of double-well potentials. Furthermore, the impacts of the time step, data size, initial values, noise ratio and different combinations of equation coefficients on the performance of EnSIL method have been explored in a systematical way. All these results confirm the efficiency and reliability of our new method for solving the inverse problems of XDEs.

Despite the successful examples presented in the main text, there are several critical issues related to the EnSIL method, which require further careful examinations.
\begin{itemize}
\item
From the perspective of information theory, since the entropy balance equation is a high compression of the original XDEs, there is no guarantee that all unknown coefficients of the XDEs could be extracted from the entropy balance equation alone. However, in most real applications, we believe the EnSIL method works.
\item
During the study, we find that the free energy is more suitable than the entropy function for solving the inverse problems of XDEs. We suspect this is probably due to the non-positiveness of the free energy dissipation rate, a manifestation of the second law of thermodynamics. However, as not all XDEs permit a free energy function, the entropy balance equation offers a universal candidate for applying our algorithm. The chemical Lorenz equation acts as a nice example.
\item
The reason for the good performance of EnSIL method in solving the inverse problems of XDEs is not fully understood, though our error analysis provides some preliminary insights into this issue. Its successful applications further suggest the original XDEs may not always be the best form for performing machine learning. Proper transformations of the equation forms may lead to more efficient and accurate learning algorithms, which have been well-known in traditional numeric calculations.
\end{itemize}



\section*
{Declaration of Competing Interest}
Authors declare that they have no conflict of interest.
\section*{Acknowledgements}
This work was supported by the National Key R\&D Program of China (Grant No. 2021YFA0719200), Guangdong Basic and Applied Basic Research Foundation (2023A1515010157), and the National Natural Science Foundation of China (21877070). L.H. thanks BIMSA for its hospitality and support. The authors would like to thank Dr. Liangrong Peng for helpful discussions.

\bibliographystyle{unsrt}
\bibliography{ref}

\end{document}